\documentstyle[11pt,amstex]{amsart}
\newcommand{\cont}{{\rm cont}}
\newcommand{\Ker}{{\rm Ker}}
\newcommand{\id}{{\rm id}}
\newcommand{\pfil}{{\rm pfil}}
\newcommand{\half}{{\rm half}}
\newcommand{\cl}{{\rm cl}}
\newcommand{\cf}{{\rm cf}}
\newcommand{\Dom}{{\rm Dom}}
\newcommand{\otp}{{\rm otp}}
\newcommand{\rest}{{\,|\grave{}\,}}
\newcommand{\lesdot}{\mathrel{\mathord{<}\!\!\raise 
0.8 pt\hbox{$\scriptstyle\circ$}}} 

\newcommand{\QED}{\hfill\vrule width 6pt height 6pt depth 0pt 
\vspace{0.1in}} 
\newcommand{\Proof}{\noindent{\sc Proof} \hspace{0.2in}} 

\newtheorem{theorem}{Theorem}[section] 
\newtheorem{claim}{Claim}[theorem]
\newtheorem{proposition}[theorem]{Proposition} 
\newtheorem{mlem}[theorem]{Main Lemma} 
\newtheorem{lemma}[theorem]{Lemma} 

\theoremstyle{definition}
\newtheorem{definition}[theorem]{Definition}

\theoremstyle{remark}
\newtheorem{remark}[theorem]{Remark}
\newtheorem{conclusion}[theorem]{Conclusion}

\title[\it The Lifting Problem With The Full Ideal]{
\vspace {3.0cm}\uppercase 
{\Large \bf The Lifting Problem With The Full Ideal}} 
\author[S. Shelah]{\uppercase {\bf S. Shelah}}
\address{Institute of Mathematics\\
The Hebrew University\\
Jerusalem 91904, Israel\\
and Rutgers University\\
Mathematics Department\\
New Brunswick, NJ 08854, USA
}
\email{shelah@@math.huji.ac.il}
\date{\today} 

\thanks{I would like to thank Alice Leonhardt for the beautiful
typing.\newline   
The research was partially supported by ``Basic Research Foundation'' of
the Israel Academy of Sciences and Humanities. Publication 636.}

\subjclass{Primary: 03E05, 28A51}

\begin{document}

\maketitle 

\bigskip
\bigskip
\bigskip
{\it Abstract.}
We show that there are a cardinal $\mu$, a $\sigma$-ideal $I\subseteq{\mathcal
P}(\mu)$ and a $\sigma$-subalgebra ${\mathcal B}$ of subsets of $\mu$
extending $I$ such that ${\mathcal B}/I$ satisfies the c.c.c. but the
quotient algebra ${\mathcal B}/I$ has no lifting.
\bigskip
\bigskip
\bigskip
\setcounter{section}{-1}

\stepcounter{section}
\subsection*{\quad 0. Introduction}
In the present paper we prove the following theorem.

\begin{theorem}
\label{0.1}
For some $\mu$ (in fact, $\mu = (2^{\aleph_0})^{++}$ suffices) there is a
$\sigma$-ideal $I$ on ${\mathcal P}(\mu)$ and a $\sigma$-subalgebra ${\frak
B}$ of ${\mathcal P}(\mu)$ extending $I$ such that ${\frak B}/I$ satisfies
the c.c.c.\ but ${\frak B}/I$ has no lifting.    
\end{theorem}
This result answers a question of David Fremlin (see chapter on measure
algebras in Fremlin \cite{Fr}). Moreover, it solves the problem of
topologizing a Category Base (see Detlefsen Szyma\'nski \cite{DS}, Morgan
\cite{Mo90}, Shilling \cite{Sl} and Szyma\'nski \cite{Sz}). 
\medskip

Note that it is well known (Mokobodzki's theorem; see Fremlin \cite{Fr}) that
under CH, if $|{\frak B}/I| \le (2^{\aleph_0})^+$ then this is impossible;
i.e.\ the quotient algebra ${\frak B}/I$ has a lifting.  

Toward the end we deal with having better $\mu$.
\medskip

I thank Andrzej Szyma\'nski for asking me the question and Max Burke and
Mariusz Rabus for corrections. 
\medskip

\noindent{\bf Notation:}\quad Our notation is rather standard. All cardinals
are assumed to be infinite and usually they are denoted by $\lambda$,
$\kappa$, $\mu$.

In Boolean algebras we use $\cap$ (and $\bigcap$), $\cup$ (and
$\bigcup$) and $-$ for the Boolean operations. 

\stepcounter{section}
\subsection*{\quad 1. The proof of Theorem \ref{0.1}}

\begin{mlem}
\label{1.1}
Suppose that
\begin{enumerate}
\item[(a)] $\mu,\lambda$ are cardinals satisfying $\mu=\mu^{\aleph_0},
\lambda\le 2^\mu$,
\item[(b)] ${\frak B}$ is a complete c.c.c.\ Boolean algebra,
\item[(c)] $x_i\in {\frak B}\setminus\{0\}$ for $i<\lambda$,
\item[(d)] for each sequence $\langle (u_i,f_i):i<\lambda\rangle$ such that
$u_i\in [\lambda]^{\le \aleph_0}$, $f_i \in {}^{\textstyle u_i}2$ there are
$n<\omega$ (but $n > 0$) and $i_0<i_1\ldots<i_{n-1}$ in $\lambda$ such that:
\begin{enumerate}
\item[$(\alpha)$] the functions $f_{i_0},\dotsc,f_{i_{n-1}}$ are compatible,
\item[$(\beta)$]  ${\frak B}\models\bigcap\limits_{\ell<n}x_{i_\ell}=0$.
\end{enumerate}
\end{enumerate}
{\em Then}
\begin{enumerate}
\item[$(\bigoplus)$] there are a $\sigma$-ideal $I$ on ${\mathcal P}(\mu)$ and
a $\sigma$-algebra ${\frak A}$ of subsets of $\mu$ extending $I$ such that
${\frak A}/I$ satisfies the c.c.c.\ and the natural homomorphism ${\frak A}
\longrightarrow {\frak A}/I$ cannot be lifted.
\end{enumerate}
\end{mlem}

\Proof Without loss of generality the algebra ${\frak B}$ has cardinality
$\lambda^{\aleph_0}$ ($\le 2^\mu$).  Let $\langle Y_b:b \in {\frak B}\rangle$ 
be a sequence of subsets of $\mu$ such that any non-trivial countable Boolean 
combination of the $Y_b$'s is non-empty (possible by \cite{EK} as $\mu =
\mu^{\aleph_0}$ and the algebra ${\frak B}$ has cardinality $\le 2^\mu$; see
background in \cite{GiSh:597}). Let ${\frak A}_0$ be the Boolean subalgebra of
${\mathcal P}(\mu)$ generated by $\{Y_b:b \in {\frak B}\}$. So $\{Y_b:b\in 
{\frak B}\}$ freely generates ${\frak A}_0$ and hence there is a unique
homomorphism $h_0$ from ${\frak A}_0$ into ${\frak B}$ satisfying $h_0(Y_b)=
b$. 

A Boolean term $\sigma$ is hereditarily countable if $\sigma$ belongs to
the closure $\Sigma$ of the set of terms $\bigcap\limits_{i<i^*} y_i$ for 
$i^*<\omega_1$ under composition and under $-y$. 
 
Let ${\mathcal E}$ be the set of all equations $\bf e$ of the form
$0=\sigma(b_0,b_1,\ldots,b_n,\ldots)_{n<\omega}$ which hold in ${\frak B}$,
where $\sigma$ is hereditarily countable.
%and holds under the substitution $y_i \mapsto x_i$
For ${\bf e}\in{\mathcal E}$ let $\cont(\bf e)$ be the set of $b\in {\frak B}$
mentioned in it (i.e.\ $\{b_n:n < \omega\}$) and let $Z_{\bold e}\subseteq\mu$
be the set $\sigma(Y_{b_0},Y_{b_1},\ldots,Y_{b_n},\ldots)_{n <\omega}$.

Let $I$ be the $\sigma$-ideal of ${\mathcal P}(\mu)$ generated by the family
$\{Z_{\bf e}:{\bf e}\in {\mathcal E}\}$ and let ${\frak A}_1$ be the Boolean
Algebra of subsets of ${\mathcal P}(\mu)$ generated by $I\cup\{Y_b:b \in
{\frak B}\}$.  

\begin{claim}
\label{star1}
$I \cap {\frak A}_0=\Ker(h_0)$.
\end{claim}

\noindent{\em Proof of the claim:}\quad Plainly $\Ker(h_0)\subseteq I\cap
{\frak A}_0$. For the converse inclusion it is enough to consider elements of 
${\frak A}_0$ of the form 
\[Y=\bigcap_{\ell=1}^{n} Y_{b_\ell}-\bigcup^{2n}_{\ell =n+1}Y_{b_\ell}.\]
If ${\frak B}\models$``$\bigcap\limits_{\ell=1}^{n} b_\ell-\bigcup\limits_{
\ell=n+1}^{2n} b_\ell=0$'' then easily $h_0(Y) = 0$. So assume that
\[{\frak B}\models\mbox{ `` }c = \bigcap_{\ell =1}^{n} b_\ell-\bigcup_{\ell =
n+1}^{2n} b_\ell \ne 0\mbox{ '',}\]
and we shall prove $Y \notin I$. Let $Z\in I$, so for some ${\bf e}_m
\in {\mathcal E}$ for $m<\omega$ we have $Z\subseteq\bigcup_{m<\omega} Z_{{\bf
e}_m}$. Let $g$ be a homomorphism from ${\frak B}$ into the 2--element Boolean
Algebra ${\frak B}_0=\{0,1\}$ such that $g(c)=1$, and $g$ respects all 
the equations ${\bf e}_m$ (including those of the form $b=\bigcup\limits_{k<
\omega} b_k$; possible by the Sikorski theorem).

By the choice of the $Y_b$'s, there is $\alpha<\mu$ such that:
\begin{quotation}
if $b\in\{b_\ell:\ell = 1,\ldots,2n\}\cup\bigcup\limits_{m<\omega}\cont({\bf
e}_m)$ then 
\[g(b)=1\Leftrightarrow\alpha\in Y_b.\]
\end{quotation}
So easily $\alpha \notin Z_{{\bf e}_m}$ for $m<\omega$, and $\alpha\in
\bigcap\limits_{\ell =1}^{n} Y_{b_\ell}\setminus\bigcup\limits_{\ell=n+1}^{2n}
Y_{b_\ell}$, so $Y$ is not a subset of $Z$. As $Z$ was an arbitrary element of
$I$ we get $Y \notin I$, so we have finished proving \ref{star1}.
\medskip

It follows from \ref{star1} that we can extend $h_0$ (the homomorphism from
${\frak A}_0$ onto ${\frak B}$) to a homomorphism $h_1$ from ${\frak A}_1$
onto ${\frak B}$ with $I=\Ker(h_1)$. Let ${\frak A}_2$ be the $\sigma$-algebra
of subsets of $\mu$ generated by ${\frak A}_1$.

\begin{claim}
\label{star2}
For every $Y \in {\frak A}_2$ there is $b \in {\frak B}$ such that $Y\equiv
Y_b \mod I$. Consequently, ${\frak A}_2={\frak A}_1$.
\end{claim}

\noindent{\em Proof of the claim:}\quad Let $Y\in {\frak A}_2$. Then $Y$ is a
(hereditarily countable) Boolean combination of some $Y_{b_\ell}$ ($\ell<
\omega$) and $Z_n$ ($n<\omega$), where $b_\ell\in {\frak B}$, $Z_n\in I$. Let
$Z_n\subseteq\bigcup\limits_{m<\omega} Z_{{\bf e}_{n,m}}$, where ${\bf e}_{n,
m}\in {\mathcal E}$, and say 
\[Y=\sigma(Y_{b_0},Z_0,Y_{b_1},Z_1,\ldots,Y_{b_n},Z_n,\ldots)_{n<\omega}.\]
Let ${\bf e}_{n,m}$ be $0=\sigma_{n,m}(b_{n,m,0},b_{n,m,1},\ldots)$. Then
clearly $\bigcup\limits_{n,m<\omega} Z_{{\bf e}_{n,m}}\in I$ (use the
definition of $I$). In ${\frak B}$, let $b=\sigma(b_0,0,b_1,0,\ldots,b_n,0,
\ldots)$ and let $\sigma^*=\sigma^*(b_0,b_1,\ldots,b_{n,m,\ell},\ldots)_{n,m,
\ell<\omega}$ be the following term  
\[\begin{array}{r}
\displaystyle\bigcup_{n,m}\sigma_{n,m}(b_{n,m,0},b_{n,m,1},\ldots)\cup(b-
\sigma(b_0,0,b_1,0,\ldots,b_m,0,\ldots))\cup\\
\cup (\sigma(b_0,0,b_1,0,\ldots,b_n,0,\ldots)-b)\cup 0.\ 
  \end{array}\]
Clearly ${\frak B}\models$``$0=\sigma^*$'', so the equation ${\bf e}$ defined
as $0=\sigma^*$ belongs to ${\mathcal E}$, and thus $Z_{{\bf e}}$ is well
defined. It follows from the definition of $\sigma^*$ that $(Y\setminus Y_b)
\cup (Y_b\setminus Y)\subseteq Z_{\bf e}\in I$.
\medskip

So we can sum up:
\begin{enumerate}
\item[(a)] $I$ is an $\aleph_1$-complete ideal of ${\mathcal P}(\mu)$,
\item[(b)] ${\frak A}_1$ is a $\sigma$-algebra of subsets of $\mu$,
\item[(c)] $I \subseteq {\frak A}_1$,
\item[(d)] $h_1$ is a homomorphism from ${\frak A}_1$ onto ${\frak B}$, with
kernel $I$,
\item[(e)] ${\frak B}$ is a complete c.c.c.\ Boolean algebra.
\end{enumerate}
This is exactly as required, so the ``only'' point left is

\begin{claim}
\label{star5}
The homomorphism $h_1$ cannot be lifted.
\end{claim}

\noindent{\em Proof of the claim:}\quad Assume that $h_1$ can be lifted, so
there is a homomorphism $g_1:{\frak B}\longrightarrow {\frak A}_1$ such that
$h_1\circ g_1 =\id_{\frak B}$.

For $i<\lambda$ let $Z_i=(g_1(x_i)-Y_{x_i})\cup (Y_{x_i}-g_1(x_i))$, so by the
assumption on $g_1$ necessarily $Z_i\in I$. Consequently we can find ${\bf
e}_{i,n}\in {\mathcal E}$ for $n<\omega$ such that $Z_i\subseteq\bigcup
\limits_{n<\omega}Z_{{\bf e}_{i,n}}$. Let $W_i=\{x_i\}\cup\bigcup\limits_{n
<\omega}\cont({\bf e}_{i,n})$, so $W_i\subseteq {\frak B}$ is countable. Let
${\frak B}'$ be the subalgebra of ${\frak B}$ generated by $\bigcup\limits_{i
< \lambda} W_i$. Clearly $|{\frak B}'|=\lambda$, so there is a one--to--one
function $t$ from $\lambda$ onto ${\frak B}'$. Put $u_i=t^{-1}(W_i)\in
[\lambda]^{\le \aleph_0}$. 

For each $i$ there is a homomorphism $f_i$ from ${\frak B}$ into the
2-element Boolean Algebra $\{0,1\}$ such that $f_i(x_i)=1$ and $f_i$ respects
all the equations ${\bf e}_{i,n}$ for $n < \omega$ (as in the proof of
\ref{star1}). Let $f'_i:u_i\longrightarrow\{0,1\}$ be defined by $f'_i(\alpha)
= f_i(t(\alpha))$. Then by clause (d) of the hypothesis there are $n<\omega$
and $i_0<\ldots<i_{n-1}<\lambda$ such that: 
\begin{enumerate}
\item[$(\alpha)$] the functions $f'_{i_0},\ldots,f'_{i_{n-1}}$ are compatible,
\item[$(\beta)$]  ${\frak B}\models$``$\bigcap\limits_{\ell<n}x_{i_\ell}=0$''.
\end{enumerate}
Hence
\begin{enumerate}
\item[$(\alpha)'$] the functions $f_{i_0}\rest W_{i_0},\ldots,f_{i_{n-1}}\rest
W_{i_{n-1}}$ are compatible\footnote{as functions, not as homomorphisms},

call their union $g$.
\end{enumerate}

Now let $\alpha<\mu$ be such that:
\begin{enumerate}
\item[$(\otimes_1)$]\quad $\ell<n\ \&\ b\in W_{i_\ell}\quad \Rightarrow\quad
[\alpha\in Y_b\Leftrightarrow g(b) = 1]$
\end{enumerate}
(it exists by the choice of the $Y_b$'s and $(\alpha)'$).

By $(\otimes_1)$ and the choice of $f_{i_\ell}$ we have:
\begin{enumerate}
\item[$(\otimes_2)$]\quad $\alpha \in Y_{x_{i_\ell}}$ 
\end{enumerate}
(because $f_{i_\ell}(x_{i_\ell})=1$) and
\begin{enumerate}
\item[$(\otimes_3)$]\quad $\alpha\notin Z_{{\bf e}_{i_\ell,n}}$ for $n<\omega$
\end{enumerate}
(because $f_{i_\ell}$ respects ${\bf e}_{i_\ell,n}$ and $\cont({\bf
e}_{i_\ell,n}) \subseteq W_{i_\ell}$) and
\begin{enumerate}
\item[$(\otimes_4)$]\quad $\alpha \notin Z_{i_\ell}$
\end{enumerate}
(by $(\otimes_3)$ as $Z_{i_\ell}\subseteq\bigcup\limits_{n<\omega}Z_{{\bf
e}_{i_\ell,n}}$).

\noindent So $\alpha\in Y_{x_{i_\ell}}\setminus Z_{i_\ell}$ and thus $\alpha
\in g_1(x_{i_\ell})$. Hence $\alpha\in\bigcap\limits_{\ell<n}g_1(x_{i_\ell})$. 
Since $g_1$ is a homomorphism we have
\[\bigcap_{\ell<n} g_1(x_{i_\ell})=g_1(\bigcap_{\ell<n} x_{i_\ell})=g_1(0)= 
\emptyset\]
(we use clause $(\beta)$ above). A contradiction. \QED$_{\ref{1.1}}$

\begin{remark}
\label{1.2}
\begin{enumerate}
\item Concerning the assumptions of \ref{1.1}, note that they seem closely
related to 
\begin{enumerate}
\item[$(\oplus_\mu)$] there is a c.c.c.\ Boolean Algebra ${\frak B}$ of
cardinality $\le \lambda$ which is not the union of $\le \mu$ ultrafilters
(i.e.\ $d({\frak B})>\mu$). 
\end{enumerate}
(See the proof of \ref{1.7} below). 
\item Concerning $(\oplus_\mu)$, by \cite{Sh:90}, if $\lambda=\mu^+$, $\mu=
\mu^{\aleph_0}$ then there is no such Boolean algebra.  By \cite{Sh:126}, it
is consistent then $\lambda=\mu^{++}\le 2^\mu$, $\aleph_0<\mu=\mu^{<\mu}$ and
$(\oplus_\mu)$ above holds using (see below) a Boolean algebra of the form
$BA(W)$, $W\subseteq [\lambda]^3$, $(\forall u_1\ne u_2\in W)(|u_1\cap u_2|\le
1)$. Hajnal, Juhasz and Szentmiklossy \cite{HaJuSz} prove the existence of a
c.c.c.\ Boolean algebra ${\frak B}$ with $d({\frak B})=\mu$ of cardinality
$2^\mu$ when there is a Jonsson algebra on $\mu$ (or $\mu$ is a limit
cardinal) using $BA(W)$, $W\subseteq [\lambda]^{<\aleph_0}$, $u\ne v\in W\quad
\Rightarrow\quad |u \cap v| < |u|/2$. The claim we need is close to this. On
the existence of Jonson cardinals (and its history) see \cite{Sh:g}. Of
course, also in \ref{1.7} if $\mu$ is not strong limit, instead ``$M$ is a
Jonsson algebra on $\mu$'' it suffices that ``$M$ is not the union of $<\mu$
subalgebras''.  Rabus Shelah \cite{RbSh:631} prove the existence of a c.c.c.\
Boolean Algebra ${\frak B}$ with $d({\frak B})=\mu$ for every $\mu$.
\end{enumerate}
\end{remark}

\begin{definition}
\label{1.3}
\begin{enumerate}
\item For a set $u$ let 
\[\begin{array}{ll}
\pfil(u)\stackrel{\rm def}{=}\{w:&w\subseteq {\mathcal P}(u),\ u\in w,\ w
\mbox{ is upward closed and}\\
\ &\mbox{if }(u_1,u_2)\mbox{ is a partition of $u$ then }u_1\in w\mbox{ or }
u_2\in w\}
  \end{array}\]
[$\pfil$ stands for ``pseudo-filter''].
\item The canonical ($\pfil$) $w$ of $u$ for a finite set $u$ is 
\[\half(u)=\{v\subseteq u:|v|\ge |u|/2\}.\] 
\item We say that $(W,{\bf w})$ is a $\lambda$-candidate if:
\begin{enumerate}
\item[(a)] $W\subseteq [\lambda]^{< \aleph_0}$,
\item[(b)] $\bf w$ is a function with domain $W$,
\item[(c)] ${\bf w}(u)\in\pfil(u)$ for $u\in W$
\item[(d)] if $v\in [\lambda]^{<\aleph_0}$ then $\cl_{(W,{\bf w})}(v)
\stackrel{\rm def}{=}\{u\in W:u \cap v\in {\bf w}(u)\}$ is finite.
\end{enumerate}
\item We say $W$ is a $\lambda$-candidate if $(W,\half\rest W$) is a
$\lambda$-candidate. 
\item Instead of $\lambda$ we can use any ordinal (or even set). 
\item We say that ${\mathcal U}\subseteq\lambda$ is $(W,{\bf w})$--closed if
for each $u\in W$
\[u\cap{\mathcal U}\in{\bf w}(u)\quad\Rightarrow\quad u\subseteq{\mathcal U}.\]
\end{enumerate}
\end{definition}

\begin{definition}
\label{1.4}
\begin{enumerate}
\item For a $\lambda$--candidate $(W,{\bf w})$ let $BA(W,{\bf w})$ be the
Boolean algebra generated by $\{x_i:i < \lambda\}$ freely except
\[\bigcap_{i\in u} x_i = 0 \qquad \mbox{ for } \qquad u \in W.\]
\item For a $\lambda$-candidate $W$, let
\[BA(W) = BA(W,\half\rest W).\]
\item For a $\lambda$-candidate $(W,{\bf w})$ let $BA^c(W,{\bf w})$ be the
completion of $BA(W,{\bf w})$; similarly $BA^c(W)$.
\end{enumerate}
\end{definition}

\begin{proposition}
\label{1.5}
Let $(W,{\bf w})$ be a $\lambda$-candidate. Then the Boolean algebra $BA(W,{
\bf w})$ satisfies the c.c.c.\ and has cardinality $\lambda$, so $BA^c(W,{\bf
w})$ satisfies the c.c.c.\ and has cardinality $\le \lambda^{\aleph_0}$. 
\end{proposition}

\Proof  Let $b_\alpha=\sigma_\alpha(x_{i_{\alpha,0}},\ldots,x_{i_{\alpha,
n_\alpha-1}})$ be nonzero members of $BA(W,{\bf w})$ (for $\alpha<\omega_1$
and $\sigma_\alpha$ a Boolean term). Without loss of generality $\sigma_\alpha
=\sigma$, $n_\alpha=n(*)$ and $i_{\alpha,0}<i_{\alpha,1}<\ldots<i_{\alpha,
n_\alpha-1}$, and $\langle\langle i_{\alpha,\ell}:\ell<n(*)\rangle:\alpha<
\omega_1\rangle$ forms a $\Delta$-system, so 
\[i_{\alpha_1,\ell_1}=i_{\alpha_2,\ell_2}\ \&\ \alpha_1\ne\alpha_2\quad
\Rightarrow\quad \ell_1=\ell_2 \ \&\ (\forall\alpha<\omega_1)(i_{\alpha,
\ell_1}=i_{\alpha_1,\ell_1}).\]
Also we can replace $b_\alpha$ by any nonzero $b'_\alpha \le b_\alpha$, so
without loss of generality for some $s_\alpha\subseteq n(*)$ ($=\{0,\ldots,
n(*)-1\}$) we have 
\[b_\alpha=\bigcap_{\ell\in s_\alpha}x_{i_{\alpha,\ell}}\cap\bigcap_{\ell\in
n(*) \setminus s_\alpha}(-x_{i_{\alpha,\ell}})>0\]
and without loss of generality $s_\alpha = s$. Put (for $\alpha<\omega_1$)
\[{\bf u}_\alpha\stackrel{\rm def}{=}\{u\in W:u\cap\{i_{\alpha,\ell}:\ell\in
s\}\in {\bf w}(u)\}\]
and note that these sets are finite (remember \ref{1.3}(3d)). Hence the sets
\[u_\alpha=\bigcup\{u:u \in {\bf u}_\alpha\}\]
are finite. Without loss of generality $\langle\{i_{\alpha,\ell}:\ell<n(*)\}
\cup u_\alpha:\alpha<\omega_1\rangle$ is a $\Delta$-system. Now let $\alpha
\ne \beta$ and assume $b_\alpha \cap b_\beta = 0$. Clearly we have
\[b_\alpha \cap b_\beta=\bigcap_{\ell\in s}(x_{i_{\alpha,\ell}}\cap
x_{i_{\beta,\ell}})\cap\bigcap_{\ell\in n(*)\setminus s} (-x_{i_{\alpha,\ell}}
\cap - x_{i_{\beta,\ell}}).\]
Note that, by the $\Delta$-system assumption, the sets $\{i_{\alpha,\ell},
i_{\beta,\ell}:\ell \in s\}$, $\{i_{\alpha,\ell},i_{\beta,\ell}:\ell\in n(*)
\setminus s\}$ are disjoint. So why is $b_\alpha\cap b_\beta$ zero? The only
possible reason is that for some $u\in W$ we have $u\subseteq\{i_{\alpha,
\ell},i_{\beta,\ell}:\ell \in s\}$.  Thus
\[u=(u\cap\{i_{\alpha,\ell}:\ell\in s\})\cup\{u\cap\{i_{\beta,\ell}:\ell\in
s\})\]
and without loss of generality $u\cap\{i_{\alpha,\ell}:\ell\in s\}\in {\bf
w}(u)$. Hence $u \in {\bf u}_\alpha$ and therefore $u\subseteq u_\alpha$. Now
we may easily finish the proof. \QED$_{\ref{1.5}}$

\begin{remark}
\label{1.6}
If we define a $(\lambda,\kappa)$--candidate weakening clause (d) to
\begin{enumerate}
\item[(d)$_\kappa$]  $v\in [\lambda]^{<\aleph_0}\quad \Rightarrow\quad
\kappa>|\{u\in W:u\cap v \in {\bf w}(u)\}|$,
\end{enumerate}
then the algebra $BA(W,{\bf w})$ satisfies the $\kappa^+$-c.c.c. 

\noindent [Why? We repeat the proof of Proposition \ref{1.5} replacing
$\aleph_1$ with $\kappa$. There is a difference only when ${\bf u}_\alpha$ has
cardinality $<\kappa$ (instead being finite) and (being the union of $<\kappa$
finite sets) also $u_\alpha$ has carinality $\mu_\alpha <\kappa$. Wlog
$\mu_\alpha =\mu<\kappa$. Clearly the set
\[S\stackrel{\rm def}{=}\{\delta<\kappa^+: \cf(\delta)=\mu^+\}\]
is a stationary subset of $\kappa^+$, so for some stationary subset $S^*$ of
$S$ and $\alpha(*)<\kappa$ we have:
\[(\forall\alpha\in S^*)\big(u_\alpha\cap\alpha\subseteq\alpha^*\quad\&\quad
u_\alpha\subseteq\min(S^*\setminus (\alpha+1))\big).\]
Let us define $u^*_\alpha=u_\alpha\cup\{i_{\alpha,\ell}:\ell\in s\}\setminus
\alpha(*)$. Wlog $\langle u^*_\alpha:\alpha\in S^*\rangle$ is a
$\Delta$--system. The rest should be clear.]
\end{remark}

\begin{theorem}
\label{1.7}
Assume that there is a Jonsson algebra on $\mu$, $\lambda = 2^\mu$, and 
\[(\forall\alpha<\mu)(|\alpha|^{\aleph_0}<\mu=\cf(\mu)).\]
Then for some $\lambda$--candidate $(W,{\bf w})$ the
Boolean algebra $BA^c(W,{\bf w})$ and $\lambda$ satisfy the assumptions {\bf
(b)}---{\bf (d)} of \ref{1.1}. 
\end{theorem}

\Proof Let $F:[\mu]^{< \aleph_0}\longrightarrow\mu$ be such that
\[(\forall A \in [\mu]^\mu)[F''([A]^{< \aleph_0}\setminus [A]^{<2})=\mu]\]
(well known and easily equivalent to the existence of a Jonsson algebra).

Let $\langle\bar{A}^\alpha:\alpha<2^\mu\rangle$ list the sequences $\bar{A}=
\langle A_i:i<\mu\rangle$ such that 
\begin{itemize}
\item \quad $A_i\in [2^\mu]^\mu$, 
\item \quad $(\forall i<\mu)(\exists\alpha)(A_i\subseteq [\mu\times\alpha,\mu
\times\alpha+\mu))$, and 
\item \quad $i<j<\mu\quad \Rightarrow\quad A_i\cap A_j=\emptyset$.
\end{itemize}
Without loss of generality we have $A^\alpha_i\subseteq\mu\times (1+\alpha)$
and each $\bar{A}$ is equal to $\bar{A}^\alpha$ for $2^\mu$ ordinals $\alpha$.
Clearly $\otp(A^\alpha_i)=\mu$. 

By induction on $\alpha<2^\mu$ we choose pairs $(W_\alpha,{\bf w}_\alpha)$ and
functions $F_\alpha$ such that 
\begin{enumerate}
\item[$(\alpha)$]  $(W_\alpha,{\bf w}_\alpha)$ is a $\mu\times(1+
\alpha)$--candidate, 
\item[$(\beta)$]   $\beta<\alpha$ implies $W_\beta=W_\alpha\cap [\mu\times(1+
\beta)]^{<\aleph_0}$ and ${\bf w}_\beta={\bf w}_\alpha\rest W_\beta$,
\item[$(\gamma)$]  $F_\alpha$ is a one-to-one function from the set 
\[\{u:u\subseteq[\mu\times (1+\alpha),\mu\times(1+\alpha +1))\mbox{ finite
with $\ge 2$ elements }\}\]
into $\bigcup\limits_{i<\mu} A^\alpha_i$,
\item[$(\delta)$]  $W_{\alpha +1} = W_\alpha\cup\{u\cup\{F_\alpha (u)\}:u\in
W^*_\alpha\}$, where 
\[W^*_\alpha=\{u:u\mbox{ a subset of }[\mu \times (1 + \alpha),\mu\times(1+
\alpha +1))\mbox{ such that }\aleph_0> |u|\ge 2\},\]
\item[$(\varepsilon)$]  for any (finite) $u\in W^*_\alpha$ we have
\[{\bf w}_{\alpha+1}(u\cup\{F_\alpha(u)\})=\{v\subseteq u\cup\{F_\alpha(u)\}:
u\subseteq v\mbox{ or } F_\alpha(u) \in v\ \&\ v\cap u\ne\emptyset\},\] 
\item[$(\zeta)$]  $F_\alpha$ is such that for any subset $X$ of $J_\alpha
= [\mu\times(1+\alpha),\mu\times(1+\alpha+1))$ of cardinality $\mu$ and $i<
\mu$ and $\gamma\in A^\alpha_i$ for some finite subset $u$ of $X$ with $\ge 2$
elements we have $F_\alpha(u)\in A^\alpha_i\setminus\gamma$.
\end{enumerate}
There is no problem to carry out the definition so that clauses
$(\beta)$--$(\zeta)$ are satisfied (to define functions $F_\alpha$ use the
function $F$ chosen at the beginning of the proof). Then $(W_\alpha,{\bf
w}_\alpha)$ is defined for each $\alpha<2^\mu$ (at limit stages $\alpha$ we
take $W_\alpha=\bigcup\limits_{\beta<\alpha} W_\beta$, ${\bf w}_\alpha=
\bigcup\limits_{\beta<\alpha} {\bf w}_\beta$, of course). 

\begin{claim}
\label{1.7.1}
For each $\alpha<2^\mu$, $(W_\alpha,{\bf w}_\alpha)$ is a $\mu\times (1+
\alpha)$--candidate.   
\end{claim}

\noindent{\em Proof of the claim:}\quad We should check the requirements of
\ref{1.3}(3). Clauses (a), (b) there are trivially satisfied. For the clause
(c) note that every element $u$ of $W_\alpha$ is of the form $u'\cup\{F_\beta
(u')\}$ for some $\beta<\alpha$ and $u'\in W_\beta^*$. Now, if $u=u_0\cup
u_1$ then either one of $u_0,u_1$ contains $u'$ or one of the two sets contains
$F_\beta(u')$ and has non-empty intersection with $u'$. In both cases we are
done. Regarding the demand (d) of \ref{1.3}(3), note that if
\[v\in [2^\mu]^{<\aleph_0},\quad u\in W_\alpha,\quad u=u'\cup \{F_\beta(
u')\},\quad u'\in W^*_\beta,\quad \beta<\alpha\]
and $v\cap u\in {\bf w}_{\beta+1}(u)$ then $v\cap u'\neq \emptyset$ and 
either\ \ $u'\subseteq v$\ \ or \ \ $F_\beta(u')\in u$.
Hence, using the fact that the functions $F_\gamma$ are one-to-one, we easily
show that for every $v\in [2^\mu]^{<\aleph_0}$ the set 
\[\{u\in W_\alpha: u\cap v\in {\bf w}_\alpha(u)\}\]
is finite (remember the definition of ${\bf w}_{\beta+1}$), finishing the
proof of the claim.
\medskip

Let $W=\bigcup\limits_{\alpha} W_\alpha$, ${\bf w}=\bigcup\limits_{\alpha}
{\bf w}_\alpha$, ${\frak B} = BA^c(W,{\bf w})$. It follows from \ref{1.7.1}
that $(W,{\bf w})$ is a $\lambda$--candidate. The main point of the proof of
the theorem is clause {\bf (d)} of the assumptions of \ref{1.1}. So let
$f_\alpha:u_\alpha\longrightarrow\{0,1\}$ for $\alpha<2^\mu$, $u_\alpha\in
[2^\mu]^{\le \aleph_0}$, be given. For each $\alpha<2^\mu$, by the assumption
that $(\forall\beta<\mu)[|\beta|^{\aleph_0}<\mu=\cf(\mu)]$ and by the
$\Delta$-lemma, we can find $X_\alpha\in[\mu]^\mu$ such that $\langle f_{\mu
\times\alpha+\zeta}:\zeta\in X_\alpha \rangle$ forms a $\Delta$-system with
heart $f^*_\alpha$. Let  
\[G=\{g:g\mbox{ is a partial function from } 2^\mu \mbox{ to }\{0,1\}\mbox{
with countable domain}\}.\]
For each $g \in G$ let $\langle\gamma(g,i):i<i(g)\rangle$ be a maximal
sequence such that $g \subseteq f^*_{\gamma(g,i)}$ and 
\[\Dom(f^*_{\gamma(g,i)})\cap\Dom(f^*_{\gamma(g,j)})=\Dom(g)\qquad\mbox{ for
$j<i$}\]
(just choose $\gamma(g,i)$ by induction on $i$).

By induction on $\zeta\le\omega_1$, we choose $Y_\zeta,G_\zeta,Z_\zeta$ and
$U_{\zeta,g}$ such that
\begin{enumerate}
\item[(a)]  $Y_\zeta\in [2^\mu]^{\le\mu}$ is increasing continuous in $\zeta$,
\item[(b)]  $Z_\zeta\stackrel{\rm def}{=} \bigcup\{\Dom(f_\gamma):(\exists
\alpha\in Y_\zeta)[\mu\times\alpha\le\gamma<\mu\times(\alpha+1)]\}$,
\item[(c)]  $G_\zeta=\{g\in G:\Dom(g)\subseteq Z_\zeta\}$,
\item[(d)]  for $g\in G_\zeta$ we have:\quad $U_{\zeta,g}$ is $\{i:i<i(g)\}$
if $i(g)<\mu^+$ and otherwise it is a subset of $i(g)$ of cardinality $\mu$
such that 
\[j\in U_{\zeta,g}\quad \Rightarrow\quad \Dom(f^*_{\gamma(g,j)})\cap Z_\zeta=
\Dom(g),\] 
\item[(e)]  $Y_{\zeta +1}=Y_\zeta\cup\{\gamma(g,j):g\in G_\zeta\mbox{ and } j
\in U_{\zeta,g}\}$. 
\end{enumerate}
Let $Y=Y_{\omega_1}$. Let $\{(g_\varepsilon,\xi_\varepsilon):\varepsilon<
\varepsilon(*)\}$, $\varepsilon(*)\leq\mu$, list the set of pairs $(g,\xi)$
such that $\xi<\omega_1$, $g\in G_\xi$ and $i(g)\ge\mu^+$. We can find
$\langle\zeta_\varepsilon:\varepsilon<\varepsilon(*)\rangle$ such that
$\langle\gamma(g_\varepsilon,\zeta_\varepsilon):\varepsilon<\varepsilon(*)
\rangle$ is without repetition and $\zeta_\varepsilon\in U_{g_\varepsilon,
\xi_\varepsilon}$. Then for some $\alpha < 2^\mu\setminus Y_{\omega_1}$ we
have  
\[(\forall\varepsilon<\varepsilon(*))(A^\alpha_\varepsilon=\{\mu\times\gamma
(g_\varepsilon,\zeta_\varepsilon)+\Upsilon:\Upsilon\in X_{\gamma(
g_\varepsilon,\zeta_\varepsilon)}\}).\] 
Now let $g=f^*_\alpha \rest Z_{\omega_1}$. Then for some $\zeta_0(*)<
\omega_1$ we have $g\in G_{\zeta_0(*)}$ and thus $U_{g,\zeta}\subseteq i(g)$
for $\zeta\in [\zeta_0(*),\omega_1)$ and $\langle\gamma(g,i):i<i(g)\rangle$
are well defined. Now, $\alpha$ exemplifies that $i(g)<\mu^+$ is impossible
(see the maximality of $i(g)$, as otherwise $i<i(g)\quad\Rightarrow\quad
\gamma(g,i)\in Y_{\zeta_0(*)+1}\subseteq Y_{\omega_1}$).

Next, for each $\gamma\in X_\alpha$, $\Dom(f_{\mu\times\alpha+\gamma})$ is
countable and hence for some $\zeta_{1,\gamma}(*)<\omega_1$ we have $\Dom(
f_{\mu\times\alpha+\gamma})\cap Z_{\omega_1}\subseteq Z_{\zeta_{1,\gamma}
(*)}$. As $\cf(\mu)>\aleph_1$ necessarily for some $\zeta_1(*)<\omega_1$ we
have that $X'_\alpha\stackrel{\rm def}{=}\{\gamma\in X_\alpha:\zeta_{1,
\gamma}(*)\le\zeta_1(*)\}\in [\mu]^\mu$, and without loss of generality
$\zeta_1(*) \ge \zeta_0(*)$. 

So for some $\varepsilon<\varepsilon(*)\le\mu$ we have $g_\varepsilon=g\ \&\
\xi_\varepsilon=\zeta_1(*)+1$. Let $\Upsilon_\varepsilon=\gamma(g_\varepsilon,
\zeta_\varepsilon)$.  Clearly
\begin{enumerate}
\item[$(*)_1$]\quad  $f^*_\alpha,f^*_{\Upsilon_\varepsilon}$ are compatible
(and countable),
\item[$(*)_2$]\quad  $\langle f_{\mu\times\alpha+\gamma}:\gamma\in X'_\alpha
\rangle$ is a $\Delta$-system with heart $f^*_\alpha$.
\end{enumerate}
So possibly shrinking $X'_\alpha$ without loss of generality
\begin{enumerate}
\item[$(*)_3$]\quad  if $\gamma\in X'_\alpha$ then $f_{\mu\times\alpha+
\gamma}$ and $f^*_{\Upsilon_\varepsilon}$ are compatible.
\end{enumerate}
For each $\gamma \in X'_\alpha$ let
\[t_\gamma=\{\beta\in X_{\Upsilon_\varepsilon}:f_{\mu\times
\Upsilon_\varepsilon+\beta}\mbox{ and } f_{\mu\times\alpha+\gamma}\mbox{ are
incompatible\/}\}.\]
As $\langle f_{\mu\times\Upsilon_\varepsilon+\beta}:\beta\in X_{
\Upsilon_\varepsilon}\rangle$ is a $\Delta$-system with heart $f^*_{
\Upsilon_\varepsilon}$ (and $(*)_3$) necessarily
\begin{enumerate}
\item[$(*)_4$]\quad $\gamma \in X'_\alpha$ implies $t_\gamma$ is countable.
\end{enumerate}
For $\gamma \in X'_\alpha$ let
\[\begin{array}{ll}
s_\gamma \stackrel{\rm def}{=}\bigcup\{u:&u\mbox{ is a finite subset of }
X'_\alpha\mbox{ and}\\
\ &F_\alpha(\{\mu\times\alpha+\beta:\beta\in u\})\mbox{ belongs to }
t_\gamma\}.
  \end{array}\] 
As $F_\alpha$ is a one-to-one function clearly
\begin{enumerate}
\item[$(*)_5$]\quad  $s_\gamma$ is a countable set.
\end{enumerate}
Hence without loss of generality (possibly shrinking $X'_\alpha$), as $\mu>
\aleph_1$,
\begin{enumerate}
\item[$(*)_6$]\quad if $\gamma_1\ne\gamma_2$ are from $X'_\alpha$ then
$\gamma_1\notin s_{\gamma_2}$. 
\end{enumerate}
By the choice of $F_\alpha$ for some finite subset $u$ of $X'_\alpha$ with at
least two elements, letting $u'\stackrel{\rm def}{=}\{\mu\times\alpha+j:j\in
u\}$ we have  
\[\beta\stackrel{\rm def}{=} F_\alpha(u')\in\{\mu\times\gamma(g_\varepsilon,
\zeta_\varepsilon)+\gamma:\gamma\in X_{\gamma(g_\varepsilon,
\zeta_\varepsilon)}\}\]
(remember $\Upsilon_\varepsilon=\gamma(g_\varepsilon,\zeta_\varepsilon))$,
so $u'\cup\{\beta\}\in W$. Thus it is enough to show that $\{f_{\mu \times
\alpha+j}:j\in u\}\cup\{f_\beta\}$ are compatible. For this it is enough to
check any two. Now, $\{f_{\mu\times\alpha+j}:j\in u\}$ are compatible as
$\langle f_{\mu\times\alpha+j}:j\in X_\alpha\rangle$ is a $\Delta$-system. So
let $j \in u$, why are $f_{\mu \times \alpha +j}$, $f_\beta$ compatible? As
otherwise $\beta-(\mu\times\Upsilon_\varepsilon)\in t_j$ and hence $u$ is a
subset of $s_j$. But $u$ has at least two elements, so there is $\gamma\in u
\setminus\{j\}$. Now $u$ is a subset of $X'_\alpha$ and this contradicts the
statement $(*)_6$ above, finishing the proof. \QED$_{\ref{1.7}}$

\begin{remark}  
In \ref{1.7}, we can also get $d(BA(W,{\bf w})) = \mu$, but this is irrelevant
to our aim. E.g.\ in this case let for $i<\mu$, $h_i$ be a partial function
from $2^\mu$ to $\{0,1\}$ such that $\Dom(h_i)\cap [\beta,\beta+\mu)$ is
finite for $\beta < 2^\mu$ and such that every finite such function is
included in some $h_i$. Choosing the $(W_\alpha,{\bf w}_\alpha)$ preserve: 
\[\{x_\beta:h_i(\beta) = 1\}\cup\{-x_\beta:h_i(\beta) = 0\}\mbox{ generates a 
filter of }BA(W_\alpha,{\bf w}_\alpha).\]
\end{remark}

\begin{conclusion}
\label{1.8}
Theorem \ref{0.1} holds.
\end{conclusion}

\Proof  By \ref{1.1}, \ref{1.7}. \QED$_{\ref{1.9}}$

\stepcounter{section}
\subsection*{\quad 2. Getting the example for $\mu=(\aleph_2)^{\aleph_0}$,
$\lambda=2^{\aleph_2}$} 
Our aim here is to show that there are $I$, ${\frak B}$ as in \ref{0.1} for
$\mu = (\aleph_2)^{\aleph_0}$. For this we shall weaken the conditions in the
Main Lemma \ref{1.1} (see \ref{1.9} below) and then show that we can get it in
a variant of \ref{1.7} (see \ref{1.10} below). More fully, by \ref{1.10} there
is a $2^{\aleph_2}$--candidate $(W,{\bf w})$ satisfying the assumptions of
\ref{1.9} except possibly clause {\bf (a)}, so $\mu$ is irrelevant in the
clauses {\bf (b)}--{\bf (f)}. Let $\mu=(\aleph_2)^{\aleph_0}=\aleph_2+
2^{\aleph_0}$ and apply \ref{1.10}. Now we get the conclusion of \ref{1.1}
as required.   

\begin{proposition}
\label{1.9}
Assume that 
\begin{enumerate}
\item[(a)]  $\mu=\mu^{\aleph_0}$, $\lambda \le 2^\mu$,
\item[(b)]  ${\frak B}$ is a complete c.c.c.\ Boolean Algebra,
\item[(c)]  $x_i\in {\frak B}\setminus\{0\}$ for $i<\lambda$, and ${\mathcal
S}\subseteq\{u\in [\lambda]^{\le\aleph_0}:(\forall i\in\lambda\setminus u)(x_i
\notin {\frak B}_u)\}$, where  ${\frak B}_u$ is the completion of $\langle
\{x_i:i \in u\}\rangle_{\frak B}$ in ${\frak B}$ (for $u\in [\lambda]^{\le
\aleph_0}$),
\item[(d)$^-$] if $i\in u_i\in [\lambda]^{\le\aleph_0}$ for $i<\lambda$, then
we can find $n<\omega$, $i_0<\ldots<i_{n-1}<\lambda$ and $u\in {\mathcal S}(
\subseteq [\lambda]^{\le \aleph_0})$ such that: 
\begin{enumerate}
\item[(i)]\ \ \ \ ${\frak B}\models$``$\bigcap\limits_{\ell<n} x_{i_\ell}=0$'',
\item[(ii)]\ \ \  $i_\ell\in u_{i_\ell}\setminus u$ for $\ell<n$,
\item[(iii)]\ \   $\langle u_{i_\ell}\setminus u:\ell<n\rangle$ are pairwise
disjoint;
\end{enumerate}
\item[(e)] $u\in{\mathcal S}\ \&\ i\in\lambda\setminus u\ \&\ y\in{\frak B}_u
\setminus\{0,1\}\quad\Rightarrow\quad y\cap x_i\ne 0\ \&\  y-x_i\ne 0$,
\item[(f)]  ${\mathcal S}$ is cofinal in $([\mu]^{<\aleph_0},{\subseteq})$

[actually, it follows from (d)$^-$].
\end{enumerate}
Then there are a $\sigma$-ideal $I$ on ${\mathcal P}(\mu)$ and a
$\sigma$-algebra ${\frak A}$ of subsets of $\mu$ extending $I$ such that
${\frak A}/I$ satisfies the c.c.c.\ and the natural homomorphism ${\frak A} 
\longrightarrow {\frak A}/I$ cannot be lifted.
\end{proposition}

\noindent{\em Remark:}\qquad  Actually we can in clause {\bf (e)} omit
``$y-x_i \ne 0$''.

\Proof  Repeat the proof of \ref{1.1} till the definition of ${\bf e}_{i,n}$
and $W_i$ in the beginning of the proof of \ref{star5} (which says that $h_2$
cannot be lifted). Then choose $u_i\in {\mathcal S}$ such that $W_i \subseteq
{\frak B}_{u_i}$ (exists by clause {\bf (f)} of our assumptions). By clause
{\bf (d)}$^-$ we can find $n<\omega$, $i_0<\ldots<i_{n-1}$ and $u\in {\mathcal
S}$ such that clauses (i),(ii),(iii) of {\bf (d)}$^-$ hold. 

\begin{claim}
For $\ell<n$, there are homomorphisms $f_{i_\ell}$ from ${\frak B}$ into $\{0,
1\}$ respecting ${\bf e}_{i_\ell,m}$ for $m<\omega$ and mapping $x_{i_\ell}$
to 1 such that $\langle f_{i_\ell}\rest (W_{i_\ell}\cap {\frak B}_u):\ell<n
\rangle$ are compatible functions.
\end{claim}

\noindent{\em Proof of the claim:}\quad E.g.\ by absoluteness it suffices to
find it in some generic extension. Let $G_u\subseteq {\frak B}_u$ be a generic
ultrafilter. Now ${\frak B}_u\lesdot {\frak B}$ and $(\forall y \in G_u)(y
\cap x_{i_\ell}>0)$ (see clause {\bf (e)}). So there is a generic ultrafilter
$G_\ell$ of ${\frak B}$ extending $G_u$ such that $x_{i_\ell}\in G_\ell$.
Define $f_{i_\ell}$ by $f_{i_\ell}(y) = 1\quad \Leftrightarrow\quad y\in
G_\ell$ for $y \in u_{i_\ell}$. By Clause (iii) of {\bf (d)}$^-$ those
functions are compatible and we finish as in \ref{1.1}.\medskip

Thus we have finished. \QED$_{\ref{1.9}}$

\begin{theorem}
\label{1.10}
In \ref{1.7} if we let e.g.\ $\mu = \aleph_2$ then we can find a
$2^\mu$-candidate $(W,{\bf w})$ such that $BA^c(W,{\bf w})$ satisfies the
clauses {\rm (b)--(f)} of \ref{1.9}. 
\end{theorem}

\Proof In short, we repeat the proof of \ref{1.7} after defining $(W,{\bf
w})$.  But now we are being given $\langle u_i:i<\lambda\rangle$, $u_i\in
[2^\mu]^{\le \aleph_0}$, $i \in u_i$. For each $\alpha<2^\mu$ (we cannot in
general find a $\Delta$-system but) we can find $u^*_\alpha$, $X_\alpha$ such
that $X_\alpha \in [\mu]^\mu,u^*_\alpha\in{\mathcal S}\subseteq [2^\mu]^{\le
\aleph_0}$ and $\langle u_{\mu \times \alpha + i} \setminus u^*_\alpha:i \in
X_\alpha \rangle$ are pairwise disjoint, and  $i\in X_\alpha\quad\Rightarrow
\quad\mu\times\alpha+i\in u_{\mu\times\alpha+i}\setminus u^*_\alpha$ and we
continue as there (replacing the functions by the sets where instead
$G_\zeta=\{g:g\in Z_\zeta,\Dom(g)\subseteq Z_\zeta\}$ we let $h_\zeta$ be a
one-to-one function from $Z_\zeta$ onto $\mu$ and $G_\zeta=\{u\subseteq
Z_\zeta:h_\zeta{}^{''}(u)\in {\mathcal S}\}$ and instead $g=f^*_\alpha\rest
Z_{\omega_1}$ let $u^*_\alpha\cap Z_{\omega_1}\subseteq Z_{\zeta_0(*)}$,
$u^*_\alpha\cap Z_{\omega_1}\subseteq v\in G_\zeta$).
\medskip

\noindent{\sc Detailed Proof} \hspace{0.2in} Let $F^*:[\mu]^{<\aleph_0}
\longrightarrow \mu$  be such that
\[(\forall A \in [\mu]^\mu)[F''([A]^{< \aleph_0}\setminus [A]^{<2})=\mu].\]
Let $\langle\bar{A}^\alpha:\alpha<2^\mu\rangle$ list the sequences $\bar{A}=
\langle A_i:i<\mu\rangle$ such that $A_i\in [2^\mu]^\mu$, $(\forall i<\mu)(
\exists\alpha)(A_i\subseteq [\mu\times\alpha,\mu\times\alpha+\mu))$ and $i<j<
\mu\quad \Rightarrow\quad A_i\cap A_j=\emptyset$. Without loss of generality
we have $A^\alpha_i\subseteq\mu\times (1+\alpha)$ and each $\bar{A}$ is equal
to $\bar{A}^\alpha$ for $2^\mu$ ordinals $\alpha$.  Clearly $\otp(A^\alpha_i)
=\mu$.

We choose by induction on $\alpha < 2^\mu$ pairs $(W_\alpha,{\bf w}_\alpha)$
and functions $F_\alpha$ such that 
\begin{enumerate}
\item[$(\alpha)$] $(W_\alpha,{\bf w}_\alpha)$ is a $\mu\times(1+
\alpha)$-candidate,
\item[$(\beta)$]  $\beta<\alpha$ implies $W_\beta=W_\alpha\cap [\mu\times (1+
\beta)]^{< \aleph_0}$, ${\bf w}_\beta={\bf w}_\alpha\rest W_\beta$,
\item[$(\gamma)$] $F_\alpha$ is a one-to-one function from 
\[\{u:u\subseteq [\mu\times (1+\alpha),\mu\times (1+\alpha+1))\mbox{ finite
with at least two elements}\}\]
into $\bigcup\limits_{i<\mu} A^\alpha_i$,
\item[$(\delta)$] $W_{\alpha+1}=W_\alpha\cup\{u\cup\{F_\alpha (u)\}:u\in
W^*_\alpha\}$, where $W^*_\alpha=\{u:u$ is a subset of $[\mu\times(1+\alpha),
\mu\times (1+\alpha +1))$ such that $\aleph_0>|u|\ge 2\}$,
\item[$(\varepsilon)$] for finite $u\in W^*_\alpha$ we have 
\[{\bf w}(u\cup \{F_\alpha(u)\})=\{v\subseteq u \cup \{F_\alpha(u)\}:u
\subseteq v \mbox{ or } F_\alpha(u)\in v\ \&\ v\cap u\ne\emptyset\},\]
\item[$(\zeta)$]  Let $F_\alpha$ be such that for any subset $X$ of $J_\alpha
= [\mu\times (1+\alpha),\mu\times(1+\alpha+1))$ of cardinality $\mu$ and $i<
\mu$ and $\gamma \in A^\alpha_i$ for some finite subset $u$ of $X$ we have
$F_\alpha(u)\in A^\alpha_i\setminus\gamma$.
\end{enumerate}
There are no difficulties in carrying out the construction and checking that
it as required. Let $W=\bigcup\limits_\alpha W_\alpha$, ${\bf w}=
\bigcup\limits_\alpha {\bf w}_\alpha$, ${\frak B} = BA^c(W,{\bf w})$. Clearly
$(W,{\bf w})$ is a $\lambda$-candidate.  

Let ${\mathcal S}^* \subseteq[\mu]^{\le \aleph_0}$ be stationary of
cardinality $\mu$.  Let 
\[{\mathcal S}'=\{u\in [\lambda]^{\le\aleph_0}:\mbox{ if }v \in W\mbox{ and } 
v\cap u\in {\bf w}(v)\mbox{ then }v\subseteq u\}.\]
Now, clause (f) holds as $(W,{\bf w})$ satisfies clause (d) of Definition
\ref{1.3}(3). As for clause (e) use Lemma \ref{1.11} below.

The main point is clause {\bf (d)}$^-$ of \ref{1.9}. So let $i\in a_i\in
[\lambda^\mu]^{\le\aleph_0}$ for $i<\lambda$ be given. For each $\alpha< 
\lambda$, as $\mu=\aleph_2$ we can find $X_\alpha\in [\mu]^\mu$ and
$a^*_\alpha \in {\mathcal S}'$ such that $\alpha\in a^*_\alpha$ and: 
\begin{enumerate}
\item[$(\otimes_\alpha)$]  $\zeta_1\ne\zeta_2\ \&\ \zeta_1\in X_\alpha\ \&\
\zeta_2\in X_\alpha\quad \Rightarrow\quad a_{\mu\times\alpha+\zeta_1}\cap
a_{\mu\times\alpha+\zeta_2}\subseteq a^*_\alpha$ and $\zeta \in X_\alpha\quad
\Rightarrow\quad \mu\times\alpha+\zeta\notin a^*_\alpha$.  
\end{enumerate}
For each $b\in [\lambda]^{\le\aleph_0}$ let $\langle\gamma(b,i):i<i(g)\rangle$
be a maximal sequence such that $\gamma(b,i)<\lambda$ and $u^*_{\gamma(b,i)}
\cap u^*_{\gamma(b,j)}\subseteq b$ and $\gamma(b,i)\notin b$ for $j<i$ (just
choose $\gamma(b,i)$ by induction on $i$).

We choose by induction on $\zeta\le\omega_1$, $Y_\zeta,h_\zeta,S_\zeta,
G_\zeta,Z_\zeta$ and $U_{\zeta,g}$ such that 
\begin{enumerate}
\item[(a)] $Y_\zeta\in [2^\mu]^{\le \mu}$ is increasing continuous in $\zeta$,
\item[(b)] $Z_\zeta$ is the minimal subset of $\lambda$ (of cardinality $\le
\mu$) which includes
\[\bigcup\{u_\gamma:(\exists\alpha\in Y_\zeta)[\mu\times\alpha\le\gamma<\mu
 \times (\alpha+1)]\}\]
and satisfies 
\[u\in W\ \&\ u\cap Z_\zeta\in {\bf w}(u)\quad\Rightarrow\quad u\subseteq
Z_\zeta,\] 
\item[(c)] $h_\zeta$ is a one-to-one function from $\mu$ onto $Z_\zeta$, and
\[G_\zeta=\{h''_\zeta(b):b\in {\mathcal S}\}\cup\bigcup\limits_{\xi<\zeta}
G_\xi.\] 
\item[(d)] for $b\in G_\zeta$ we have $U_{\zeta,b}$ is $\{i:i<i(b)\}$ if $i(b)
<\mu^+$ and otherwise is a subset of $i(b)$ of cardinality $\mu$ such that 
\[j\in U_{\zeta,b}\quad \Rightarrow \quad\Dom(f^*_{\gamma(b,j)})\cap Z_\zeta
\subseteq b,\]
\item[(e)]  $Y_{\zeta+1}=Y_\zeta\cup\{\gamma(b,j):b\in G_\zeta\mbox{ and } j
\in U_{\zeta,b}\}$. 
\end{enumerate}
Again, there is no problem to carry out the definition (e.g.\ $|Z_\zeta|\le
\mu$ by clause (d) of \ref{1.3}(3)). Let $Y=Y_{\omega_1}$. Let $\{(
b_\varepsilon,\xi_\varepsilon):\varepsilon<\varepsilon(*)\le\mu\}$ list the
set of pairs $(b,\xi)$ such that $\xi<\omega_1$, $b\in G_\xi$ and $i(b)\ge
\mu^+$. We can find $\langle\zeta_\varepsilon:\varepsilon<\varepsilon(*)
\rangle$ such that $\langle\gamma(b_\varepsilon,\zeta_\varepsilon):\varepsilon
<\varepsilon(*)\rangle$ is without repetition and $\zeta_\varepsilon\in
U_{b_\varepsilon,\xi_\varepsilon}$, $\varepsilon(*)\leq\mu$. So for some
$\alpha<2^\mu\setminus Y_{\omega_1}$ we have 
\[(\forall\varepsilon<\varepsilon(*))(A^\alpha_\varepsilon=\{\mu\times \gamma(
b_\varepsilon,\zeta_\varepsilon)+\Upsilon:\Upsilon\in X_{\gamma(b_\varepsilon,
\zeta_\varepsilon)}\}.\]
Now, let $b_0=a^*_\alpha\cap Z_{\omega_1}$, so for some $\zeta_0(*)<\omega_1$
we have $b_0\subseteq Z_{\zeta_0(*)}$. As $a^*_\alpha$ is countable and
$G_\zeta \subseteq [Z_\zeta]^{\le \aleph_0}$ is stationary (and the closure
property of $Z_\zeta$) there is $b^*\in {\mathcal S}'$ such that $b\stackrel{
\rm def}{=} b^*\cap Z_{\zeta_0(*)}$ belongs to $G_\zeta$ and $a^*_\alpha
\subseteq b^*$ and so $U_{b,\zeta}\subseteq i(b)$ for $\zeta\in [\zeta_0(*), 
\omega_1)$ and $\langle\gamma(b,i):i<i(b)\rangle$ are well defined. Now
$\alpha$ exemplified $i(b)<\mu^+$ is impossible (see the maximality as
otherwise $i<i(b)\quad \Rightarrow\quad\gamma(b,i)\in Z_{\zeta_0(*)+1}
\subseteq Z_{\omega_1}$). 

As for each $\gamma\in X_\alpha$, the set $a_{\mu\times\alpha+\gamma}$ is
countable,  for some $\zeta_{1,\gamma}(*)<\omega_1$ we have $a_{\mu\times
\alpha+\gamma}\cap Z_{\omega_1}\subseteq Z_{\zeta_{1,\gamma}(*)}$.  Since
$\cf(\mu)>\aleph_1$ necessarily for some $\zeta_1(*)<\omega_1$ we have 
\[X'_\alpha\stackrel{\rm def}{=}\{\gamma\in X_\alpha:\zeta_{1,\gamma}(*)\le
\zeta_1(*)\}\in [\mu]^\mu\]
and without loss of generality $\zeta_1(*)\ge\zeta_0(*)$. Thus for some
$\varepsilon < \mu$ we have $b_\varepsilon=b\ \&\ \xi_\varepsilon=\zeta_1(*)+
1$. Let $\Upsilon_\varepsilon=\gamma(b_\varepsilon,\zeta_\varepsilon)$. 
Clearly 
\begin{enumerate}
\item[$(*)_1$] $a^*_\alpha,a^*_{\Upsilon_\varepsilon}$ are countable,
\item[$(*)_2$] $\gamma\in X'_\alpha\quad \Rightarrow\quad \mu\times\alpha+
\gamma\in a_\gamma$, 
\item[$(*)_3$] $\gamma_1\ne\gamma_2\ \&\ \gamma_1\in X'_\alpha\ \&\ \gamma_2
\in X'_\alpha\quad \Rightarrow\quad a_{\mu\times\alpha+\gamma_1}\cap a_{\mu
  \times\alpha+\gamma_2}\subseteq b^*$. 
\end{enumerate}
So possibly shrinking $X'_\alpha$ without loss of generality
\begin{enumerate}
\item[$(*)_4$] if $\gamma\in X'_\alpha$ then $a_{(\mu\times\alpha+\gamma)}
\cap a^*_{\Upsilon_\varepsilon}\subseteq b^*$. 
\end{enumerate}
For each $\gamma \in X'_\alpha$ let
\[t_\gamma=\{\beta\in X_{\Upsilon_\varepsilon}:a_{(\mu\times
\Upsilon_\varepsilon+\beta)}\cap a_{(\mu\alpha+\gamma)}\not\subseteq
b^*\}.\]
As $\langle f_{(\mu\times\Upsilon_\varepsilon+\beta)}:\beta\in X_{
\Upsilon_\varepsilon}\rangle$ was chosen to satisfy
$(\otimes_{\Upsilon_\varepsilon})$ (and $(*)_3$) necessarily 
\begin{enumerate}
\item[$(*)_5$] $\gamma\in X'_\alpha$ implies $t_\gamma$ is countable.
\end{enumerate}
For $\gamma \in X'_\alpha$ let
\[\begin{array}{ll}
s_\gamma\stackrel{\rm def}{=}\bigcup\{u:&u\mbox{ is a finite subset of }
X'_\alpha \mbox{ and}\\
\ &F_\alpha(\{\mu\times\alpha+\beta:\beta\in u\})\mbox{ belongs to }t_\gamma
\}.\end{array}\]
As $F_\alpha$ is a one-to-one function clearly
\begin{enumerate}
\item[$(*)_6$] $s_\gamma$ is a countable set.
\end{enumerate}
So without loss of generality (possibly shrinking $X'_\alpha$ using $\mu>
\aleph_1$) 
\begin{enumerate}
\item[$(*)_7$] if $\gamma_1\ne\gamma_2$ are from $X'_\alpha$ then $\gamma_1
\notin s_{\gamma_2}$. 
\end{enumerate}
By the choice of $F_\alpha$, for some finite subset $u$ of $X'_\alpha$ with at
least two elements, letting $u'\stackrel{\rm def}{=} \{\mu\times\alpha+j:j\in
u\}$ we have 
\[\beta\stackrel{\rm def}{=} F_\alpha(u')\in \{\mu\times\gamma(b_\varepsilon,
\zeta_\varepsilon)+\gamma:\gamma\in X_{\gamma(b_\varepsilon,
\zeta_\varepsilon)}\}.\]
Hence $u'\cup\{\beta\}\in W$, so it is enough to show that $\{a_{\mu\times
\alpha+j}:j \in u\}\cup\{a_\beta\}$ are pairwise disjoint outside $b^*$. For
the first it is enough to check any two. Now, $\{f_{\mu\times\alpha +j}:j
\in u\}$ are O.K.\ by the choice of $\langle f_{\mu\times\alpha + j}:j\in
X_\alpha\rangle$. So let $j\in u$. Now, $a_{\mu \times \alpha +j}$, $a_\beta$
are O.K., otherwise $\beta-(\mu\times\Upsilon_\varepsilon)\in t_j$ and hence
$u$ is a subset of $s_j$ but $u$ has at least two elements and is a subset of
$X'_\alpha$ and this contradicts the statement $(*)_6$ above and so we are
done. \QED$_{\ref{1.10}}$

\begin{lemma}
\label{1.11}
Let $(W,{\bf w})$ be a $\lambda$-candidate. Assume that $u\subseteq\lambda$
and $u=\cl_{(W,{\bf w})}(u)$ (see Definition \ref{1.3}(1),(d)) and let
$W^{[u]}=W \cap [u]^{<\aleph_0}$ and ${\bf w}^{[u]}={\bf w}\rest
W^{[u]}$. Furthermore suppose that $(W,{\bf w})$ is non-trivial (which holds
in all the cases we construct), i.e.
\begin{enumerate}
\item[$(*)$]\qquad  $i\in v\in W\quad \Rightarrow\quad v\setminus\{i\}\in {\bf
w}(v)$. 
\end{enumerate}
Then:
\begin{enumerate}
\item $(W^{[u]},{\bf w}^{[u]})$ is a $\lambda$-candidate (here $u=\cl_{(W,{\bf
w})}(u)$ is irrelevant);
\item $BA(W^{[u]},{\bf w}^{[u]})$ is a subalgebra of $BA(W,{\bf w})$, moreover
$BA(W^{[u]},{\bf w}^{[u]})\lesdot BA(W,{\bf w})$;
\item if $i\in\lambda\setminus u$ and $y\in BA(W^{[u]},{\bf w}^{[u]})$ then
\[y\ne 0\quad \Rightarrow\quad y\cap x_i>0\ \&\ y-x_i>0;\]
\item $BA^c(W^{[u]},{\bf w}^{[u]})\lesdot BA^c(W,{\bf w})$.
\end{enumerate}
\end{lemma}

\Proof  1)\quad Trivial.
\medskip

\noindent 2)\quad {\em The first phrase}:\quad if $f_0$ is a homomorphism from
$BA(W^{[u]},{\bf w}^{[u]})$ to the Boolean Algebra $\{0,1\}$ we define a
function $f$ from $\{x_\alpha:\alpha<\lambda\}$ to $\{0,1\}$ by $f(x_\alpha)$
is $f_0(x_\alpha)$ if $\alpha\in u$ and is zero otherwise. Now 
\[v\in W\quad\Rightarrow\quad (\exists\alpha \in v)(f(x_\alpha) = 0).\]
Why?  If $v \subseteq u$, then $v \in W^{[u]}$ and ``$f_0$ is a homomorphism",
so $f_0(\bigcap\limits_{\alpha\in v} x_\alpha)=0$. Hence $(\exists \alpha\in
v)(f_0(x_\alpha)=0)$ and hence $(\exists\alpha \in v)(f(x_\alpha)=0)$. If $v
\not\subseteq u$, then choose $\alpha\in v\setminus u$, so $f(x_\alpha)=0$.

So $f$ respects all the equations involved in the definition of $BA(W,{\bf
w})$ hence can be extended to a homomorphism $\hat f$ from $BA(W,{\bf w})$ to
$\{0,1\}$. Easily $f_0 \subseteq \hat f$ and so we are done.
\medskip

As for {\em the second phrase}, let $z\in BA(W,{\bf w})$, $z>0$ and we shall
find $y\in BA(W^{[u]},{\bf w}^{[u]})$, $y>0$ such that
\[(\forall x)[x\in BA(W^{[u]},{\bf w}^{[u]})\ \&\ 0<x\le y\quad \Rightarrow
\quad x\cap z\ne 0).\]
We can find disjoint finite subsets $s_0,s_1$ of $\lambda$ such that $0<z'\le
z$ where $z'=\bigcap\limits_{\alpha\in s_1} x_\alpha\cap\bigcap\limits_{\alpha
\in s_0} (-x_\alpha)$.  Let 
\[t=\bigcup\{v:v\in W\mbox{ a finite subset of }\lambda\mbox{ and }v\cap s_0
\in {\bf w}(v)\} \cup s_0 \cup s_1.\]
We know that $t$ is finite. We can find a partition $t_0,t_1$ of $t$ (so $t_0
\cap t_1=\emptyset$, $t_0\cup t_1=t$) such that $s_0\subseteq t_0$ and $s_1
\subseteq t_1$ and $y^*=\bigcap\limits_{\alpha\in t_1} x_\alpha\cap
\bigcap\limits_{\alpha\in t_0}(-x_\alpha)>0$. Note that $y\stackrel{\rm def}{=}
\bigcap\limits_{\alpha\in u\cap t_1} x_\alpha\cap \bigcap\limits_{\alpha\in u
\cap t_0}(-x_\alpha)$ is $> 0$ and, of course, $y\in BA(W^{[u]},{\bf
w}^{[u]})$. We shall show that $y$ is as required. So assume $0<x\le y$, $x\in
BA(W^{[u]},{\bf w}^{[u]})$. As we can shrink $x$, without loss of generality,
for some disjoint finite $r_0,r_1\subseteq u$ we have $t\cap u\subseteq r_0
\cup r_1$ and $x=\bigcap\limits_{\alpha\in r_1} x_\alpha\cap\bigcap\limits_{
\alpha\in r_0}(-x_\alpha)$, so clearly $t_1\cap u\subseteq r_1$, $t_0\cap u
\subseteq r_0$. 

We need to show $x\cap z\ne 0$, and for this it is enough to show that $x\cap
z' \ne 0$. Now, it is enough to find a function $f:\{x_\alpha:\alpha<\lambda\}
\longrightarrow\{0,1\}$ respecting all the equations in the definition of
$BA(W,{\bf w})$ such that $\hat{f}$ maps $x\cap z'$ to $1$. So let
$f(x_\alpha)=1$ for $\alpha\in r_1\cup s_1$ and $f(x_\alpha)=0$ otherwise. If
this is O.K., fine as $f\rest r_0$, $f \rest s_0$ are identically zero and
$f\rest r_1$, $f\rest s_1$ are identically one. If this fails, then for some
$v\in {\bf w}$ we have $v\subseteq r_1\cup s_1$. But then $v\cap r_1\in {\bf
w}(v)$ or $v\cap s_1\in {\bf w}(v)$. Now if $v \cap r_1\in {\bf w}(w)$ as $r_1
\subseteq u$ necessarily $v \subseteq u$, but $v\subseteq r_1\cup s_1$ and
$s_1 \cap u\subseteq t_1\subseteq r_1$, so $v\subseteq r_1$ is a contradiction
to $x>0$. Lastly, if $v\cap s_1\in {\bf w}(v)$, then $v\subseteq t$ so as
$v\subseteq r_1\cup s_1$ we have $v\subseteq s_1\cup (t\cap r_1)$ and so
$v\subseteq s_1\cup t_1$ and hence $v\subseteq t_1$ --- a contradiction to
$y^*>0$. So $f$ is O.K.\ and we are done. 
\medskip

\noindent 3)\quad Let $f_0$ be a homomorphism from $BA(W^{[u]},{\bf w}^{[u]})$
to the trivial Boolean Algebra $\{0,1\}$. For $t\in\{0,1\}$ we define a
function $f$ from $\{x_\alpha:\alpha < \lambda\}$ to $\{0,1\}$ by
\[f(x_\alpha)=\left\{\begin{array}{lll}
f_0(x_\alpha)&\mbox{if }&\alpha\in u\\
t            &\mbox{if }&\alpha=i\\
0            &\mbox{if }&\alpha\in\lambda\setminus u\setminus\{i\}.
                     \end{array}\right.\]
Now $f$ respects the equations in the definition of $BA(W,{\bf w})$. Why? Let
$v\in W$. We should prove that $(\exists\alpha \in v)(f(\alpha)=0)$. If $v
\subseteq u$, then 
\[f\rest\{x_\alpha:\alpha\in v\}=f_0\rest\{x_\alpha:\alpha \in v\}\quad\mbox{
and}\]
\[0=f_0(0_{BA(W^{[u]},{\bf w}^{[u]})})=f_0(\bigcap_{\alpha\in v}x_\alpha)=
\bigcap_{\alpha \in v} f_0(x_\alpha),\]
so $(\exists\alpha \in v)(f_0(x_\alpha)=0)$. If $v\not\subseteq u\cup\{i\}$
let $\alpha\in v\setminus u\setminus\{i\}$, so $f(x_\alpha)=0$ as required.

So we are left with the case $v\subseteq u\cup \{i\}$, $v\not\subseteq
u$. Then by the assumption $(*)$, $v\cap u=v\setminus\{i\}\in{\bf w}(v)$ and
$v \subseteq u$, a contradiction. 
\medskip
 
\noindent 4)\quad Follows. \QED$_{\ref{1.11}}$

\begin{remark}
\label{1.12}
We can replace $\aleph_0$ by say $\kappa =\cf(\kappa)$ (so in \ref{1.10}, $\mu
=\kappa^{++}$, in \ref{1.7}, $(\forall\alpha<\mu)(|\alpha|^{<\kappa}<\mu=
\cf(\mu))$.
\end{remark}

\end{document}